\documentclass{article}%
\usepackage{amsmath}
\usepackage{amsfonts}
\usepackage{amssymb}
\usepackage{graphicx}%
\providecommand{\U}[1]{\protect\rule{.1in}{.1in}}
\newtheorem{theorem}{Theorem}

\newtheorem{corollary}[theorem]{Corollary}

\newtheorem{example}[theorem]{Example}

\newtheorem{lemma}[theorem]{Lemma}

\newtheorem{proposition}[theorem]{Proposition}
\newtheorem{remark}[theorem]{Remark}

\begin{document}

\begin{center}
\bigskip{\Large Central limit theorem for Fourier transform and periodogram of
random fields}

\bigskip Magda Peligrad and Na Zhang

\bigskip
\end{center}

Department of Mathematical Sciences, University of Cincinnati, PO Box 210025,
Cincinnati, Oh 45221-0025, USA. \texttt{ }

email: peligrm@ucmail.uc.edu

email: zhangn4@mail.uc.edu

\bigskip

\textbf{Abstract.} In this paper we show that the limiting distribution of the
real and the imaginary part of the Fourier transform of a stationary random
field is almost surely an independent vector with Gaussian marginal
distributions, whose variance is, up to a constant, the field's spectral
density. The dependence structure of the random field is general and we do not
impose any restrictions on the speed of convergence to zero of the
covariances, or smoothness of the spectral density. The only condition
required is that the variables are adapted to a commuting filtration and are
regular in some sense. The results go beyond the Bernoulli fields and apply to
both short range and long range dependence. They can be easily applied to
derive the asymptotic behavior of the periodogram associated to the random
field. The method of proof is based on new probabilistic methods involving
martingale approximations and also on borrowed and new tools from harmonic
analysis. Several examples to linear, Volterra and Gaussian random fields will
be presented.

\bigskip


\bigskip


\section{Introduction}

The discrete Fourier transform, defined as%
\begin{equation}
S_{n}(t)=\sum_{k=1}^{n}\mathrm{e}^{\mathrm{i}kt}X_{k}\,, \label{Four}%
\end{equation}
where $\mathrm{i}=\sqrt{-1}$ is the imaginary unit, plays an essential role in
the study of stationary time series $(X_{j})_{j\in Z}$ of centered random
variables with finite second moment, adapted to a filtration $(\mathcal{F}%
_{u})_{u\in Z}$.

The periodogram, introduced as a tool by Schuster (1898), is essential for the
estimation of the spectral density of the stationary processes. It is defined
by
\begin{equation}
I_{n}(t)=\frac{1}{2\pi n}|S_{n}(t)|^{2},\ t\in\lbrack-\pi,\pi). \label{Per}%
\end{equation}
There is a vast literature concerning these statistics. They are often used to
determine hidden periodicities. Denote by $\lambda$ the Lebesgue measure on
the real line. In Peligrad and Wu (2010), it was proved a surprising result,
that, under ergodicity and a very mild regularity condition, for $\lambda
-$almost all frequencies $t$, the random variables $\operatorname{Re}%
S_{n}(t)/\sqrt{n}$ and$\,\operatorname{Im}S_{n}(t)/\sqrt{n}$ are
asymptotically independent identically distributed random variables with
normal distribution, mean $0$ and variance $\pi f(t)$. Here $f$ is the
spectral density of $(X_{j})_{j\in Z}$. The regularity condition, namely
$E(X_{0}|\mathcal{F}_{-\infty})=0$ $a.s.$, is a mild restriction of dependence
and accommodates large classes of sequences with short and long range
dependence. This result implies that for $\lambda-$almost all $t,$ the
periodogram $I_{n}(t)$ converges in distribution to $f(t)\chi^{2},$ where
$\chi^{2}$ has\textrm{ }a chi-square distribution with $2$ degrees of freedom,
even in the case of processes with long memory. The proof of this result is
based on the celebrated Carleson theorem (1966) on almost sure convergence of
Fourier transforms, combined with martingale approximations and Fourier analysis.

In this paper we analyze the asymptotic properties of the Fourier transform
for random fields. Let $d$ be a positive integer. We start with a strictly
stationary random field $(X_{\mathbf{u}})_{\mathbf{u}\in Z^{d}}$ of square
integrable and centered random variables. We introduce the discrete Fourier
transform for random fields by the rotated sum%
\[
S_{\mathbf{n}}(\mathbf{t)=}%
{\displaystyle\sum\limits_{\mathbf{1}\leq\mathbf{u}\leq\mathbf{n}\ }}
\mathrm{e}^{\mathrm{i}\mathbf{u}\cdot\mathbf{t}}X_{\mathbf{u}},
\]
where we have $1\leq\mathbf{u\leq n}$ and $\mathbf{t}\in I=[-\pi,\pi)^{d}.$ By
$\mathbf{u\leq n}$ we understand $\mathbf{u=}(u_{1},...,u_{d})$,
$\mathbf{n=}(n_{1},...,n_{d})$ and $1\leq u_{1}\mathbf{\leq}n_{1},$..., $1\leq
u_{d}\mathbf{\leq}n_{d}\mathbf{.}$ Also $\mathbf{u}\cdot\mathbf{t=}u_{1}%
t_{1}+...+u_{d}t_{d}.$

For a weakly stationary random field the spectral analysis was initiated in
several papers by Helson and Lowdenslager (1962), Kallianpur et al. (1990) and
Francos et al. (1995). These papers stress the huge difficulties when one
tries to extend the results from sequences of random variables to random
fields. One of the difficulty is that for random fields the future and past do
not have a unique interpretation. Also, many of the important spectral
analysis results relevant to the proofs, do not fully extend to double indexed
sequences, including the celebrated Fej\'{e}r-Lebesgue Theorem (cf. Bary,
1964, p. 139) or the Carleson theorem (1966) (see Fefferman, (1971,b)). To
compensate for the lack of ordering of the filtration we utilize the notion of
commuting filtration. Such filtrations have a certain Markovian character. For
instance, for $d=2,$ we can start with a stationary random field with
independent rows or columns which generate a commuting filtration. Then, we
construct a stationary random field which is a function of the initial one.

The main result of the paper is a natural extension from sequences of random
variables, indexed by integers, to random fields of the result of Peligrad and
Wu (2010). Under certain regularity conditions we shall prove that, almost
surely in $\mathbf{t}\in I,$ both the real and imaginary part of
$S_{\mathbf{n}}(\mathbf{t)/}\sqrt{n_{1}...n_{d}}$ converge to independent
normal variables whose variance is, up to a multiplicative constant, the
spectral density of the random field, denoted by $f(\mathbf{t})$. The
ergodicity condition is imposed to only one of the directions of the random field.

The periodogram, has the following extension to the random field:
\[
I_{\mathbf{n}}(\mathbf{t})={\frac{1}{({2\pi)}^{d}{n}_{1}...n_{d}}}\left\vert
{\displaystyle\sum\limits_{\mathbf{1}\leq\mathbf{u}\leq\mathbf{n}\ }}
\mathrm{e}^{\mathrm{i}\mathbf{u}\cdot\mathbf{t}}X_{\mathbf{u}}\right\vert
^{2},\quad\mathbf{t}\in I.
\]
Our result is that, for almost all frequencies $\mathbf{t}\in I$, the limiting
distribution of $I_{\mathbf{n}}(\mathbf{t})$ is $f(\mathbf{t})\chi^{2}(2),$
where $\chi^{2}(2)$ is a chi-square distribution with two degrees of freedom.

The proof is based on a new, interesting representation for the spectral
density in terms of projection operators, which is the most important tool for
establishing our result. The proof also involves a martingale approximation
for random fields as well as laws of large numbers for Fourier sums, which
have interest in themselves.

We consider two types of summations. The first result is for summations of the
variables in a multi-dimensional cube. The reason we first restrict ourselves
to summations indexed by the cubes is due to the relation between our results
and optimal results available in harmonic analysis. For example, for $d=2,$
Theorem 1 in Marcinkiewicz and Zygmund (1939) shows that the
Fej\'{e}r-Lebesgue theorem holds for spectral densities in $L_{1}$ when the
summation is taken over rectangles of size $m\times n,$ provided that
$m,n,\rightarrow\infty$ such that $m/n\leq a$ and $n/m\leq a$ for some
positive number $a$. This result fails when the summation is taken over
general rectangles. However, if the summation is taken over the sets $1\leq
u_{1}\mathbf{\leq}n\mathbf{,}$ $1\leq u_{2}\mathbf{\leq}m\mathbf{\ }$\ where
$n\geq m\rightarrow\infty,$ one should assume the integrability of
$f(\mathbf{u})\ln^{+}f(\mathbf{u})$ as a minimal condition for the validity of
the Fej\'{e}r-Lebesgue theorem (see Jessen et al. 1935). We shall also give a
result in this context, where the summation is taken over unrestricted rectangles.

When dealing with random fields the notation can become rather complicated.
This is the reason why, for presenting the material, we implemented the
following strategy: We treat first the case $d=2.$ Then, we mention the small
differences for treating the general case of multi-dimensional index set by
using the mathematical induction.

Our paper is organized as follows. In a preliminary section we review several
facts about limiting variance of the Fourier series, introduce the notions of
stationary random fields and commuting filtrations. In Section 3 we obtain a
representation of the spectral density in terms of projection operators, which
extends a recent result by Lifshitz and Peligrad (2015) beyond the setting of
Bernoulli shifts. We also state and prove our main results on the limiting
distribution of double indexed, random Fourier sums. The extension to general
index set is given in Section 4. Section 5 is dedicated to examples, such as
functions of Gaussian sequences, linear and nonlinear random fields with
independent innovations. It is remarkable that the only condition required for
the validity of our results for linear or Volterra random fields with
independent innovations is equivalent to merely the existence of these fields.
In a supplementary section we prove two laws of large numbers and other lemmas
about commuting filtrations.

Our paper joins the recent increased interest in finding martingale methods
for random fields, initiated by Rosenblatt (1972) and continued by Gordin
(2009). We would like to mention several remarkable papers in this direction.
For instance, the paper by Voln\'{y} and Wang (2014) treated projection
conditions and orthomartingales. Voln\'{y} (2015) discovered that the fields
of stationary orthomartingales require the ergodicity in only one of the
directions of the field as a necessary condition for the CLT. Cuny et al.
(2015) treated dynamical system via projection conditions. The paper of El
Machkouri et al. (2013) deals with random fields which are functions of i.i.d.
Wang and Woodroofe (2013), Peligrad and Zhang (2017) and Giraudo
(2017)\ treated the Maxwell-Woodroofe condition. Also, in the context of
dynamical systems, the CLT\ for Fourier transform for random fields was stated
in Cohen and Conze (2013) for $K$-systems. It should also be mentioned that a
central limit theorem for periodogram of random fields was obtained by Miller
(1995)\ under mixing conditions. All these papers were inspirational.

\section{Preliminaries.}

\textbf{Spectral density and limiting variance}

\bigskip

We call the complex valued zero mean field of random variables $(X_{\mathbf{m}%
})_{\mathbf{m}\in Z^{2}}$ defined on a probability space $(\Omega
,\mathcal{K},P),$ weakly stationary (or second order stationary), if there are
complex numbers $\gamma(\mathbf{m}),$ $\mathbf{m}\in Z^{2}$, such that for all
$\mathbf{u},\mathbf{v}\in Z^{2},$%
\[
\mathrm{cov}(X_{\mathbf{u}},X_{\mathbf{v}})=E(X_{\mathbf{u}}\bar
{X}_{\mathbf{v}})=\gamma(\mathbf{u}-\mathbf{v}).
\]
In the context of weakly stationary random fields it is known that there
exists a unique measure on $I=[-\pi,\pi)^{2}$, such that
\[
\gamma(\mathbf{u)}=\int_{I}\mathrm{e}^{\mathrm{i}\mathbf{u}\cdot\mathbf{x}%
}F(\mathrm{d}\mathbf{x}),\quad\text{for all}\,\,\mathbf{u}\in Z^{2}\,,
\]
where $\mathbf{u}\cdot\mathbf{x}$ is the inner product. If $F$ is absolutely
continuous with respect to Lebesgue measure $\lambda^{2}$ on $I=[-\pi,\pi
)^{2}$ then, the Radon-Nikodym derivative $f$ of $F$ with respect to the
Lebesgue measure is called \textit{spectral density} ($F(\mathrm{d}%
\mathbf{t})=f(\mathbf{t})\mathrm{d}\mathbf{t}$)$,$ and we have%
\[
\gamma(\mathbf{u)}=\int_{I}\mathrm{e}^{\mathrm{i}\mathbf{u}\cdot\mathbf{x}%
}f(\mathbf{x})\mathrm{d}\mathbf{x},\quad\text{for all}\,\,\mathbf{u}\in
Z^{2}\,.
\]
The variance of partial sums on rectangles is%

\[
E|S_{\mathbf{n}}(\mathbf{t)|}^{2}=%
{\displaystyle\sum\limits_{1\leq\mathbf{u},\mathbf{v\leq n}}}
\gamma(\mathbf{u}-\mathbf{v)}e^{\mathrm{i}\mathbf{t}\cdot(\mathbf{u}%
-\mathbf{v})}.
\]
Well-known computations show that%
\begin{align*}
E|S_{\mathbf{n}}(\mathbf{t)}|^{2}  &  =%
{\displaystyle\sum\limits_{\mathbf{1}\leq\mathbf{u},\mathbf{v\leq n}}}
e^{\mathrm{i}\mathbf{t}\cdot(\mathbf{u}-\mathbf{v})}\int_{I}e^{\mathrm{i}%
\mathbf{x}\cdot(\mathbf{u-v)}}f(\mathbf{x})\mathrm{d}\mathbf{x}\\
&  =\int_{I}%
{\displaystyle\sum\limits_{\mathbf{1}\leq\mathbf{u},\mathbf{v\leq n}}}
e^{\mathrm{i}\mathbf{x}\cdot(\mathbf{u}-\mathbf{v})}f(\mathbf{x}%
-\mathbf{t})\mathrm{d}\mathbf{x.}%
\end{align*}
So, with the notation $\mathbf{x}=(x_{1},x_{2}),$ one can rewrite
\[
\frac{1}{n_{1}n_{2}}E|S_{\mathbf{n}}(\mathbf{t)}|^{2}=\int_{I}K_{n_{1}}%
(x_{1})K_{n_{2}}(x_{2})f(\mathbf{x}-\mathbf{t})\mathrm{d}\mathbf{x,}%
\]
where $K_{n}(x)$ is the Fej\'{e}r Kernel%
\[
K_{n}(x)=\sum\nolimits_{|j|<n}(1-\frac{|j|}{n})\mathrm{e}^{\mathrm{i}jx}.
\]
Furthermore, by Theorem 1 in Marcinkiewicz and Zygmund (1939), for
$\lambda^{2}-$almost all $\mathbf{t}$ in $I,$ we obtain a limiting
representation for the spectral density, namely
\begin{equation}
\lim_{n\rightarrow\infty}\frac{1}{n^{2}}E|S_{n,n}(\mathbf{t)}|^{2}=(2\pi
)^{2}f(\mathbf{t).} \label{limvar}%
\end{equation}
If in addition $f(\mathbf{u})\ln^{+}f(\mathbf{u})$ is integrable, then%
\begin{equation}
\lim_{n_{1}\geq n_{2}\rightarrow\infty}\frac{1}{n_{1}n_{2}}E|S_{n_{1},n_{2}%
}(\mathbf{t)}|^{2}=(2\pi)^{2}f(\mathbf{t).} \label{limvarln}%
\end{equation}

\bigskip

\textbf{Stationary random fields and stationary filtrations}

\bigskip

In order to construct stationary filtrations, we shall start with a strictly
stationary real valued random field $\mathbf{\xi}=(\xi_{\mathbf{u}%
})_{\mathbf{u}\in Z^{2}}$, defined on a probability space $(\Omega
,\mathcal{K},P)$ and\ define the filtrations%
\begin{equation}
\mathcal{F}_{k,\ell}=\sigma(\xi_{j,u}:j\leq k,\text{ }u\leq\ell
).\label{def fitration}%
\end{equation}
To ease the notation, sometimes the conditional expectation will be denoted by%
\[
E_{a,b}X=E(X|\mathcal{F}_{a,b}).
\]
We shall consider that the filtration is commuting in the sense that
\begin{equation}
E_{u,v}E_{a,b}X=E_{a\wedge u,b\wedge v}X.\mathbb{\ }\label{pcf}%
\end{equation}
It is remarkable that, by Problem 34.11 in Billingsley (1995) (see Lemma
\ref{commuting}), condition (\ref{pcf})\ is equivalent to the apparently
weaker condition: for $a\geq u$ and $X$ integrable we have$\ $%
\begin{equation}
E_{u,v}E_{a,b}X=E_{u,b\wedge v}X.\label{tcf}%
\end{equation}
Now we introduce the stationary random field $(X_{\mathbf{m}})_{\mathbf{m}\in
Z^{2}}$, in the following way. We define first%
\[
X_{\mathbf{0}}=g((\mathbf{\xi}_{\mathbf{u}})_{\mathbf{u}\in Z^{2}}).
\]
where $g:R^{Z^{2}}\rightarrow\mathbb{C}$ and $\mathbf{0}=(0,0)$.

Without restricting the generality we shall define $(\mathbf{\xi}_{\mathbf{u}%
})_{\mathbf{u}\in Z^{2}}$ in a canonical way on the probability space $\Omega$
$=R^{Z^{2}}$, endowed with the $\sigma-$field, $\mathcal{B}(\Omega),$
generated by cylinders. Then, if $\omega=(x_{\mathbf{v}})_{\mathbf{v}\in
Z^{2}}$ we define $\mathbf{\xi}_{\mathbf{u}}^{\prime}(\omega)=x_{\mathbf{u}}$.
We construct a probability measure $P^{\prime}$ on $\mathcal{B}(\Omega)$ such
that for all $B\in\mathcal{B}(\Omega)$ and any $m$ and $\mathbf{u}%
_{1},...,\mathbf{u}_{m}$ we have%
\[
P^{\prime}((x_{\mathbf{u}_{1}},...,x_{\mathbf{u}_{m}})\in B)=P((\mathbf{\xi
}_{\mathbf{u}_{1}},...,\mathbf{\xi}_{\mathbf{u}_{m}})\in B).
\]
The new sequence $(\mathbf{\xi}_{\mathbf{u}}^{\prime})_{\mathbf{u}\in Z^{2}}$
is distributed as $(\mathbf{\xi}_{\mathbf{u}})_{\mathbf{u}\in Z^{2}}$ and
re-denoted $(\mathbf{\xi}_{\mathbf{u}})_{\mathbf{u}\in Z^{2}}$. We shall also
re-denote $P^{\prime}$ as $P.$ Now on $R^{Z^{2}}$ we introduce the operators%
\[
T^{\mathbf{u}}((x_{\mathbf{v}})_{\mathbf{v}\in Z^{2}})=(x_{\mathbf{v+u}%
})_{\mathbf{v}\in Z^{2}}.
\]
Two of them will play an important role in our paper namely, when
$\mathbf{u=}(1,0)$ and when $\mathbf{u=}(0,1).$ By interpreting the indexes as
notations for the lines and columns of a matrix, we shall call%

\[
T((x_{u,v})_{(u,v)\in Z^{2}})=(x_{u+1,v})_{(u,v)\in Z^{2}}%
\]
the vertical shift and%
\[
S((x_{u,v})_{(u,v)\in Z^{2}})=(x_{u,v+1})_{(u,v)\in Z^{2}}%
\]
the horizontal shift.\ Then define
\begin{equation}
X_{j,k}=g(T^{j}S^{k}(\mathbf{\xi}_{\mathbf{u}})_{\mathbf{u}\in Z^{2}}).
\label{defXfield}%
\end{equation}

\section{Results and proofs}

\textbf{Spectral density representation in terms of projections}

\bigskip

In this section we first find a useful representation of the spectral density
for regular functions and commuting filtrations. It extends a result of
Lifshitz and Peligrad (2015) beyond the case of Bernoulli shifts. The proof
follows the same lines as in Lifshitz and Peligrad (2015). We shall point out
the differences and give it here for completeness, clarification and
equivalent definitions.

For an integrable random variable $X$, we introduce the projection operators
by defining%
\[
P_{\tilde{0},0}X=(E_{0,0}-E_{-1,0})X
\]%
\[
P_{0,\tilde{0}}X=(E_{0,0}-E_{0,-1})X.
\]
Note that, by (\ref{pcf}), we have%
\[
{\mathcal{P}}_{{\mathbf{0}}}X=P_{\tilde{0},0}P_{0,\tilde{0}}X=P_{0,\tilde{0}%
}P_{\tilde{0},0}X.
\]
and by stationarity, for all $u,v\in Z$
\begin{equation}
{\mathcal{P}}_{u,v}X=E_{u,v}X-E_{u,v-1}X-E_{u-1,v}X+E_{u-1,v-1}X. \label{proj}%
\end{equation}

Define ${\mathcal{F}}_{-\infty,m}=\cap_{u\in Z}{\mathcal{F}}_{u,m}$ and
${\mathcal{F}}_{m,-\infty}=\cap_{v\in Z}{\mathcal{F}}_{m,v}.$

Now let $X_{\mathbf{0}}$ be defined as before, in $L_{2}(\Omega,\mathcal{F}%
,${$P$}$).$ Then we have the following orthogonal representation
\[
X_{\mathbf{0}}={\sum\limits_{{\mathbf{u}}\in{\mathbb{J}}_{n,m}}}{\mathcal{P}%
}_{{\mathbf{u}}}X_{\mathbf{0}}+R_{nm}+U_{nm},
\]
where $n$ and $m$ are two positive integers, ${\mathbb{J}}_{n,m}%
=[-n,...,n]\times\lbrack-m,...,-m]$,
\[
R_{nm}=E(X_{\mathbf{0}}|{\mathcal{F}}_{-n-1,m})+E(X_{\mathbf{0}}|{\mathcal{F}%
}_{n,-m-1})-E(X_{\mathbf{0}}|{\mathcal{F}}_{-n-1,-m-1}),
\]
and%
\[
U_{nm}=X_{\mathbf{0}}-E(X_{\mathbf{0}}|{\mathcal{F}}_{n,m}).
\]
We assume the following two regularity conditions
\begin{equation}
E(X_{\mathbf{0}}|{\mathcal{F}}_{-\infty,0})=0\text{ }a.s.\text{ and
}E(X_{\mathbf{0}}|{\mathcal{F}}_{-\infty,0})\text{ }a.s. \label{reg}%
\end{equation}
Note that%
\[
E(X_{\mathbf{0}}|{\mathcal{F}}_{-n-1,-m-1})=E(E(X_{\mathbf{0}}|{\mathcal{F}%
}_{0,-m-1})|{\mathcal{F}}_{-n-1,0}).
\]
By passing to the limit and using the reverse martingale theorem and arguments
similar to Theorem 34.2 (V) in Billingsley (1995), we obtain that%
\[
\lim_{n\rightarrow\infty}\lim_{m\rightarrow\infty}R_{nm}=0\text{ \ }a.s.\text{
and in }L_{2}.
\]
Since $X_{\mathbf{0}}$ is measurable with respect to $\vee_{\mathbf{u}\in
Z^{2}}\mathcal{F}_{\mathbf{u}}=\mathcal{F}_{\infty,\infty},$ by the martingale
convergence theorem,%
\[
\lim_{n\rightarrow\infty}\lim_{m\rightarrow\infty}U_{nm}=0\text{ \ }a.s.\text{
and in }L_{2}.
\]
Therefore%
\[
X_{\mathbf{0}}=\lim_{n\rightarrow\infty}\lim_{m\rightarrow\infty}%
\sum\limits_{j=-n}^{n}\sum\limits_{k=-m}^{m}{\mathcal{P}}_{-j,-k}%
X_{\mathbf{0}}\text{ \ }a.s.\text{ and in }L_{2}.
\]
We shall denote this limit by%
\begin{equation}
X_{\mathbf{0}}=\sum_{\mathbf{u}\in Z^{2}}{\mathcal{P}}_{{\mathbf{u}}%
}X_{{\mathbf{0}}}. \label{Xrep}%
\end{equation}
Note that for $\mathbf{u}\neq\mathbf{v}$ and for all $X$ and $Y$ in
$L_{2}(\Omega,{\mathcal{K}},${$P$}$)$ we have%

\begin{equation}
{\mathrm{cov}}({\mathcal{P}}_{\mathbf{u}}X,{\mathcal{P}}_{\mathbf{v}}Y)=0.
\label{ortho}%
\end{equation}
Observe also that, by taking into account (\ref{Xrep}), (\ref{ortho}) and
stationarity, we have
\begin{equation}
\sum_{\mathbf{u}\in Z^{2}}E|{\mathcal{P}}_{{\mathbf{u}}}X_{{\mathbf{0}}}%
|^{2}=E|X_{{\mathbf{0}}}|^{2}<\infty. \label{s}%
\end{equation}

We would like now to define a random variable which will be used to
characterize the spectral density of random fields. For random variables this
was achieved in Peligrad and Wu (2010) by using Carleson (1966) and also Hunt
and Young (1974) theorems. For random fields these theorems do not hold in
general. We could use instead a weaker form of them or, as an alternative, an
iterated procedure.

In the sequel we shall use the notation $||X||^{2}=E|X|^{2}.$

\bigskip

\textbf{Martingale difference construction}

\bigskip

We start from the identity (\ref{s}) and note that this identity implies%
\begin{equation}
\sum_{\mathbf{u}\in Z^{2}}|{\mathcal{P}}_{{\mathbf{u}}}X_{{\mathbf{0}}}%
|^{2}<\infty\text{ \ }P-a.s. \label{star 2}%
\end{equation}
Let $\Omega^{\prime}\subset\Omega$ with $P(\Omega^{\prime})=1$ be such that
the convergence above holds for all $\omega\in\Omega^{\prime}$. By the main
theorem in Fefferman (1971 a) for convergence of double Fourier series, we
obtain the almost sure convergence in the following sense: \ For $\omega
\in\Omega^{\prime}$ we have
\[
\sum_{\mathbf{u}\in Z^{2}}\mathrm{e}^{-{\mathrm{i}\mathbf{j}}\cdot{\mathbf{t}%
}}{\mathcal{P}}_{{\mathbf{0}}}X_{{\mathbf{j}}}=\lim_{n\rightarrow\infty}%
\sum_{{\mathbf{j}}\in{\mathbb{I}}_{n}}\mathrm{e}^{-{\mathrm{i}}\,{\mathbf{j}%
}\cdot{\mathbf{t}}}{\mathcal{P}}_{{\mathbf{0}}}X_{{\mathbf{j}}}\text{
\ \ }\lambda^{2}-a.e.\text{,}%
\]
where ${\mathbb{I}}_{n}=[-n,n]\times\lbrack-n,n].$

By Fubini Theorem, for almost all ${\mathbf{t\in}}I,$ we also have that%
\begin{equation}
\sum_{\mathbf{u}\in Z^{2}}\mathrm{e}^{-{\mathrm{i}}\,{\mathbf{j}}%
\cdot{\mathbf{t}}}{\mathcal{P}}_{{\mathbf{0}}}X_{{\mathbf{j}}}=\lim
_{n\rightarrow\infty}\sum_{{\mathbf{j}}\in{\mathbb{I}}_{n}}\mathrm{e}%
^{-{\mathrm{i}}\,{\mathbf{j}}\cdot{\mathbf{t}}}{\mathcal{P}}_{{\mathbf{0}}%
}X_{{\mathbf{j}}}\text{ \ \ }P-a.s.\text{ } \label{star1}%
\end{equation}
Furthermore, by relation (1) in Fefferman's paper (1971 a) and by
(\ref{star 2}), for a positive constant $C,$ we have that
\[
\int_{I}\sup_{n}|\sum_{{\mathbf{j}}\in{\mathbb{I}}_{n}}\mathrm{e}%
^{-{\mathrm{i}}\,{\mathbf{j}}\cdot{\mathbf{t}}}{\mathcal{P}}_{{\mathbf{0}}%
}X_{{\mathbf{j}}}|^{2}d\mathbf{t\leq}C\sum_{\mathbf{u}\in Z^{2}}|{\mathcal{P}%
}_{{\mathbf{u}}}X_{{\mathbf{0}}}|^{2}\text{ \ }P-a.s.
\]
Whence, by integrating and using (\ref{s}), we obtain \
\[
E\int_{I}\sup_{n}|\sum_{{\mathbf{j}}\in{\mathbb{I}}_{n}}\mathrm{e}%
^{-{\mathrm{i}}\,{\mathbf{j}}\cdot{\mathbf{t}}}{\mathcal{P}}_{{\mathbf{0}}%
}X_{{\mathbf{j}}}|^{2}d\mathbf{t\leq}CE\sum_{\mathbf{u}\in Z^{2}}%
|{\mathcal{P}}_{{\mathbf{u}}}X_{{\mathbf{0}}}|^{2}\leq C||X_{\mathbf{0}}%
||^{2}.
\]
It follows that for almost all ${\mathbf{t\in}}I$
\[
E(\sup_{n}|\sum_{{\mathbf{j}}\in{\mathbb{I}}_{n}}\mathrm{e}^{-{\mathrm{i}%
}\,{\mathbf{j}}\cdot{\mathbf{t}}}{\mathcal{P}}_{{\mathbf{0}}}X_{{\mathbf{j}}%
}|^{2})<\infty.
\]
By the dominated convergence theorem, the convergence in (\ref{star1})\ also
holds in $L_{2}.$

Let us denote by
\begin{equation}
D_{{\mathbf{0}}}({\mathbf{t)=}}\lim_{n\rightarrow\infty}\sum_{{\mathbf{j}}%
\in{\mathbb{I}}_{n}}\mathrm{e}^{-{\mathrm{i}}\,{\mathbf{j}}\cdot{\mathbf{t}}%
}{\mathcal{P}}_{{\mathbf{0}}}X_{{\mathbf{j}}}\text{ \ }P-a.s.\text{ and in
}L_{2}. \label{defd0}%
\end{equation}

In the next theorem we point out a representations for the spectral density by
using definition (\ref{defd0}).

\begin{theorem}
\label{specden}Let $\left(  X_{{\mathbf{k}}}\right)  _{{\mathbf{k}}\in{Z}^{2}%
}$ be a stationary sequence defined by (\ref{defXfield}) and the filtration
$(\mathcal{F}_{{\mathbf{k}}})_{{\mathbf{k}}\in{Z}^{2}}$ is commuting as in
(\ref{pcf}). Assume that the second moment is finite and the regularity
condition (\ref{reg}) is satisfied. Then, the sequence $(X_{{\mathbf{k}}%
})_{{\mathbf{k}}\in{Z}^{2}}$ has spectral density which has the
representation:%
\begin{equation}
f({\mathbf{t}})=\frac{1}{(2\pi)^{2}}E\,|D_{{\mathbf{0}}}({\mathbf{t)}}%
|^{2}\qquad{\mathbf{t}}\in I. \label{rep1}%
\end{equation}

\end{theorem}

\textbf{Proof of Theorem \ref{specden}.} Let us compute the covariance of
$X_{{\mathbf{k}}}$ and $X_{{\mathbf{0}}}$. By using the projection
decomposition in (\ref{Xrep}),\ written for both $X_{{\mathbf{k}}}$ and
$X_{{\mathbf{0}}}$, together with the orthogonality of the projections in
(\ref{ortho}) and stationarity, we have for all ${\mathbf{k}}\in${$Z$}$^{2}$,
\begin{gather}
{\mathrm{cov}}(X_{{\mathbf{k}}},X_{{\mathbf{0}}})={\mathrm{cov}}%
({\sum\limits_{{\mathbf{j}}\in Z^{2}}}{\mathcal{P}}_{{\mathbf{j}}%
}X_{{\mathbf{k}}},{\sum\limits_{\mathbf{u}\in Z^{2}}}{\mathcal{P}%
}_{{\mathbf{u}}}X_{{\mathbf{0}}})\label{cov1}\\
={\sum\limits_{\mathbf{j}\in Z^{2}}}{\mathrm{cov}}({\mathcal{P}}_{{\mathbf{j}%
}}X_{{\mathbf{k}}},{\mathcal{P}}_{{\mathbf{j}}}X_{{\mathbf{0}}})={\sum
\limits_{\mathbf{j}\in Z^{2}}}{\mathrm{cov}}({\mathcal{P}}_{\mathbf{0}%
}X_{{\mathbf{k}}+{\mathbf{j}}},{\mathcal{P}}_{\mathbf{0}}X_{{\mathbf{j}}%
}).\nonumber
\end{gather}
Let us analyze the function $f({\mathbf{t}})$ defined in (\ref{rep1}). By
Fubini theorem and (\ref{s}) we have
\[
\int_{I}f({\mathbf{t}})\mathrm{d}\mathbf{t}=\frac{1}{(2\pi)^{2}}E\int_{I}%
|\sum_{\mathbf{j}\in Z^{2}}\mathrm{e}^{-{\mathrm{i}}\,{\mathbf{j}}%
\cdot{\mathbf{t}}}{\mathcal{P}}_{{\mathbf{0}}}X_{{\mathbf{j}}}|^{2}%
\mathrm{d}\mathbf{t}=E\sum_{\mathbf{j}\in Z^{2}}|{\mathcal{P}}_{{\mathbf{j}}%
}X_{{\mathbf{0}}}|^{2}<\infty.
\]
Now, let us compute the Fourier coefficients of $f({\mathbf{t}})$. For every
${\mathbf{k}}\in${$Z$}$^{2},$ by the definition of $f({\mathbf{t}})$ and
Fubini theorem we have
\[
\int_{I}{\mathrm{e}}^{{\mathrm{i}}\,{\mathbf{k}}\cdot{\mathbf{t}}%
}f({\mathbf{t}})\,\mathrm{d}\mathbf{t}=\frac{1}{(2\pi)^{2}}\sum_{\mathbf{u}%
,{\mathbf{j}}\in Z^{2}}\,E\int_{I}{\mathrm{e}}^{{\mathrm{i}}\,({\mathbf{k}%
}-{\mathbf{j+u}})\cdot{\mathbf{t}}}{\mathcal{P}}_{{\mathbf{0}}}X_{{\mathbf{j}%
}}\overline{{\mathcal{P}}_{{\mathbf{0}}}X_{{\mathbf{u}}}}\mathrm{d}%
\mathbf{t}.
\]
By using the orthogonality of the exponential functions, we obtain
\begin{align*}
\int_{I}{\mathrm{e}}^{{\mathrm{i}}\,{\mathbf{k}}\cdot{\mathbf{t}}%
}f({\mathbf{t}})\mathrm{d}\mathbf{t}  &  =\sum_{{\mathbf{j}},{\mathbf{u}}\in
Z^{2}}\mathrm{cov}\left(  {\mathcal{P}}_{{\mathbf{0}}}X_{{\mathbf{j}}%
},{\mathcal{P}}_{\mathbf{0}}X_{{\mathbf{u}}}\right)  {\mathbf{1}%
}_{\{{\mathbf{k}}-{\mathbf{j+u=0}}\}}\\
&  =\sum_{{\mathbf{u}}\in Z^{2}}\mathrm{cov}\left(  {\mathcal{P}}%
_{{\mathbf{0}}}X_{{\mathbf{u+k}}},{\mathcal{P}}_{\mathbf{0}}X_{{\mathbf{u}}%
}\right)  .
\end{align*}
Now, comparing this expression with (\ref{cov1}) we see that $f$ in formula
(\ref{rep1})\ is the spectral density for $(X_{{\mathbf{k}}})_{{\mathbf{k}}%
\in{Z}^{2}}$. \ \ $\square$ \bigskip

\begin{remark}
For defining the spectral density iterated limits are also possible. By
applying Carleson (1966) and Hunt and Young (1972) theorems twice,
consecutively in each variable, one can show that the following limits exist:
for $\lambda^{2}-$almost all $\mathbf{t}\in\lbrack-\pi,\pi)^{2},$ we can
define a random variable $\tilde{D}_{{\mathbf{0}}}({\mathbf{t)}}$ in the
following sense%
\[
\tilde{D}_{{\mathbf{0}}}({\mathbf{t)=}}\lim_{n\rightarrow\infty}%
\lim_{m\rightarrow\infty}\sum\limits_{u_{1}=-n}^{n}\sum\limits_{u_{2}=-m}%
^{m}{\mathcal{P}}_{{\mathbf{0}}}(X_{u_{1},u_{2}})\mathrm{e}^{-{\mathrm{i}%
}\,{\mathbf{u}}\cdot{\mathbf{t}}}\text{ }P\text{-a.s. and in }L_{2}.
\]
Similarly, we can also define the other iterated limit
\[
\hat{D}_{{\mathbf{0}}}({\mathbf{t)}}=\lim_{m\rightarrow\infty}\lim
_{n\rightarrow\infty}\sum\limits_{u_{1}=-n}^{n}\sum\limits_{u_{2}=-m}%
^{m}{\mathcal{P}}_{{\mathbf{0}}}(X_{u_{1},u_{2}})\mathrm{e}^{-{\mathrm{i}%
}\,{\mathbf{u}}\cdot{\mathbf{t}}}\text{ }P-a.s.\text{ and in }L_{2}.
\]
Also, we can obtain the following alternative definitions for the spectral
density:
\[
f({\mathbf{t}})=\frac{1}{(2\pi)^{2}}E\,|\hat{D}_{{\mathbf{0}}}({\mathbf{t)}%
}|^{2}=\frac{1}{(2\pi)^{2}}E\,|\tilde{D}_{{\mathbf{0}}}({\mathbf{t)}}|^{2}.
\]
Note that in all the characterizations of $f({\mathbf{t}})$ the limits commute
with the integrals.
\end{remark}

\begin{remark}
\label{rem}By (\ref{pcf}) and its definition, $D_{{\mathbf{0}}}({\mathbf{t)}}$
is a martingale difference in each coordinate%
\[
E_{0,-1}D_{{\mathbf{0}}}({\mathbf{t)=0}}\text{ and }E_{-1,0}D_{{\mathbf{0}}%
}({\mathbf{t)=0}}\text{ \ }P-a.s.
\]

\end{remark}

We are ready to state our main result. Everywhere in the paper $\Rightarrow$
denotes convergence in distribution.

\begin{theorem}
\label{CLT}Assume that $\left(  X_{{\mathbf{k}}}\right)  _{{\mathbf{k}}\in
{Z}^{2}}$ and $(\mathcal{F}_{{\mathbf{k}}})_{{\mathbf{k}}\in{Z}^{2}}$ are as
in Theorem \ref{specden}. In addition, assume that one of the shifts $T$ or
$S$ is ergodic. Then, for $\lambda^{2}-$almost all $\mathbf{t\in}I,$
\[
\frac{1}{n}(\operatorname{Re}S_{n,n}(\mathbf{t}),\operatorname{Im}%
S_{n,n}(\mathbf{t}))\Rightarrow(N_{1},N_{2})\text{ as }n\rightarrow\infty,
\]
where $N_{1},N_{2}$ are i.i.d. normally distributed random variables with mean
$0$ and variance $2\pi^{2}f(\mathbf{t}),$ where $f(\mathbf{t})$ is the
spectral density of the sequence $\left(  X_{{\mathbf{k}}}\right)
_{{\mathbf{k}}\in{Z}^{2}}$. Furthermore, if $f(\mathbf{u})\ln^{+}%
f(\mathbf{u})$ is integrable then
\[
\frac{1}{\sqrt{n_{1}n_{2}}}(\operatorname{Re}S_{n_{1},n_{2}}(\mathbf{t}%
),\operatorname{Im}S_{n_{1},n_{2}}(\mathbf{t}))\Rightarrow(N_{1},N_{2})\text{
as }n_{1}\wedge n_{2}\rightarrow\infty
\]
with $N_{1},N_{2}$ as above.
\end{theorem}

\textbf{Proof of Theorem \ref{CLT}. }This proof has several steps. Let us
point out the idea of the proof. First we show that the proof can be reduced
to random variables with continuous spectral density. Then, we construct a
random field which is a martingale difference in each coordinate and has the
same limiting distribution as the original sigma field. To validate this
approximation we shall use the limiting variance given in (\ref{limvar}) and
(\ref{limvarln})\ along with the representation of the spectral density given
in Theorem \ref{specden}. The result will follow by obtaining the central
limit theorem for the martingale random field. To fix the ideas let us assume
that the shift $S$ is ergodic.

\bigskip

\textbf{Martingale approximation}

\bigskip

Let us recall the definition of $D_{\mathbf{0}}(\mathbf{t})$ given in
(\ref{defd0})\ and introduce a new notation:
\begin{equation}
D_{\mathbf{0}}^{(\ell)}(\mathbf{t})=\sum_{{\mathbf{j}}\in\mathbb{I}_{\ell}%
}{\mathcal{P}}_{{\mathbf{0}}}(X_{{\mathbf{j}}}){\mathrm{e}}^{-{\mathrm{i}%
}\,{\mathbf{j}}\cdot{\mathbf{t}}}\rightarrow\sum_{{\mathbf{j}}\in Z^{2}%
}{\mathcal{P}}_{{\mathbf{0}}}(X_{{\mathbf{j}}}){\mathrm{e}}^{-{\mathrm{i}%
}\,{\mathbf{j}}\cdot{\mathbf{t}}}=D_{\mathbf{0}}(\mathbf{t})\text{
}P-a.s.\text{ and in }L_{2}. \label{limit}%
\end{equation}

Note that $D_{\mathbf{0}}(\mathbf{t})$ and $D_{\mathbf{0}}^{(\ell)}%
(\mathbf{t})$ are functions of $(\mathbf{\xi}_{\mathbf{u}})_{\mathbf{u}\in
Z^{2}}.$ By using stationarity and translation operators $T$ and $S$ we define
$D_{\mathbf{k}}^{(\ell)}(\mathbf{t})$ and $D_{\mathbf{k}}(\mathbf{t})$ for any
$\mathbf{k}\mathbb{\in}${$Z$}$^{2}.\ $Note that, by Remark \ref{rem}, both
$(D_{u,v}^{(\ell)}(\mathbf{t}))$ and $(D_{u,v}(\mathbf{t}))$ are
coordinate-wise martingale differences with respect to the filtrations
$(\mathcal{F}_{\infty,v})_{v}$ and $(\mathcal{F}_{u,\infty})_{u}$
respectively$.$

For almost all $\mathbf{t\in}I$ we shall approximate $S_{\mathbf{n}%
}(\mathbf{t)}$ by the martingale%
\begin{equation}
M_{\mathbf{n}}(\mathbf{t)=}\sum_{\mathbf{j}=\mathbf{1}}^{\mathbf{n}%
}{\mathrm{e}}^{{\mathrm{i}}\,{\mathbf{j}}\cdot{\mathbf{t}}}D_{\mathbf{j}%
}(\mathbf{t}). \label{def M}%
\end{equation}
To validate this approximation, we first consider the situation when
$\mathbf{n}=(n,n).$ Define the martingale%
\begin{equation}
M_{\mathbf{n}}^{(\ell)}(\mathbf{t)=}\sum_{\mathbf{j}=\mathbf{1}}^{\mathbf{n}%
}{\mathrm{e}}^{{\mathrm{i}}\,{\mathbf{j}}\cdot{\mathbf{t}}}D_{\mathbf{j}%
}^{(\ell)}(\mathbf{t}) \label{def ML}%
\end{equation}
and, for $\mathbf{t}^{\prime}$ fixed, the "proper"\ Fourier series in
$\mathbf{t,}$%

\[
M_{\mathbf{n}}^{(\ell)}(\mathbf{t},\mathbf{t}^{\prime}\mathbf{)=}%
\sum_{\mathbf{j}=\mathbf{1}}^{\mathbf{n}}{\mathrm{e}}^{{\mathrm{i}%
}\,{\mathbf{j}}\cdot{\mathbf{t}}}D_{\mathbf{j}}^{(\ell)}(\mathbf{t}^{\prime
}).
\]
Note that we can bound%
\[
|S_{\mathbf{n}}(\mathbf{t)}-M_{\mathbf{n}}(\mathbf{t)|}^{2}\mathbf{\leq
}3\mathbf{(}|S_{\mathbf{n}}(\mathbf{t)}-M_{\mathbf{n}}^{(\ell)}(\mathbf{t}%
,\mathbf{t}^{\prime}\mathbf{)|}^{2}\mathbf{+|}M_{\mathbf{n}}^{(\ell
)}(\mathbf{t},\mathbf{t}^{\prime}\mathbf{)-}M_{\mathbf{n}}^{(\ell
)}(\mathbf{t)|}^{2}\mathbf{+|}M_{\mathbf{n}}^{(\ell)}(\mathbf{t)-}%
M_{\mathbf{n}}(\mathbf{t)|}^{2}\mathbf{).}%
\]
By (\ref{limvar}), for almost all $\mathbf{t}{\mathbf{\in}}I$ \
\[
\lim_{n\rightarrow\infty}\frac{1}{n^{2}}E|S_{\mathbf{n}}(\mathbf{t)}%
-M_{\mathbf{n}}^{(\ell)}(\mathbf{t,t}^{\prime}\mathbf{)|}^{2}=(2\pi
)^{2}f^{(\ell)}({\mathbf{t,t}}^{\prime}),
\]
where $f^{(\ell)}({\mathbf{t,t}}^{\prime})$ is the spectral density of
$(X_{\mathbf{k}}-D_{\mathbf{k}}^{(\ell)}(\mathbf{t}^{\prime}))_{k}.$

By using the representation (\ref{rep1}) given in Theorem \ref{specden}, and
taking into account that ${\mathcal{P}}_{{\mathbf{0}}}D_{\mathbf{j}}^{(\ell
)}(\mathbf{t}^{\prime})=0$ $P-a.s.$ for ${\mathbf{j}}\in Z^{2}$ with
${\mathbf{j}}\neq\mathbf{0}$ we obtain
\[
f^{(\ell)}({\mathbf{t,t}}^{\prime})=\frac{1}{(2\pi)^{2}}{E}\,|\sum
_{{\mathbf{j}}\in Z_{2}}{\mathcal{P}}_{{\mathbf{0}}}X_{{\mathbf{j}}}%
\mathrm{e}^{-{\mathrm{i}}\,{\mathbf{j}}\cdot{\mathbf{t}}}-D_{\mathbf{0}%
}^{(\ell)}(\mathbf{t}^{\prime})|^{2}=\frac{1}{(2\pi)^{2}}{E}\,|D_{\mathbf{0}%
}(\mathbf{t})-D_{\mathbf{0}}^{(\ell)}(\mathbf{t}^{\prime})|^{2}.
\]
On the other hand, by the orthogonality of the projections,%
\begin{equation}
\frac{1}{n^{2}}E\mathbf{|}M_{\mathbf{n}}^{(\ell)}(\mathbf{t},\mathbf{t}%
^{\prime}\mathbf{)-}M_{\mathbf{n}}^{(\ell)}(\mathbf{t)|}^{2}=E\mathbf{|}%
D\mathbf{_{\mathbf{0}}^{(\ell)}(\mathbf{t}^{\prime})-}D\mathbf{_{\mathbf{0}%
}^{(\ell)}(\mathbf{t})|}^{2} \label{two}%
\end{equation}
and%
\[
\frac{1}{n^{2}}E\mathbf{|}M_{\mathbf{n}}^{(\ell)}(\mathbf{t)-}M_{\mathbf{n}%
}(\mathbf{t)|}^{2}=E\mathbf{|}D_{\mathbf{0}}^{(\ell)}(\mathbf{t)-}%
D_{\mathbf{0}}(\mathbf{t)|}^{2}.
\]
So, overall, by the above considerations,
\begin{gather*}
\lim\sup_{n\rightarrow\infty}\frac{1}{n^{2}}E|S_{\mathbf{n}}(\mathbf{t)}%
-M_{\mathbf{n}}(\mathbf{t)|}^{2}\leq3({E}\,|D_{0}(\mathbf{t})-D_{\mathbf{0}%
}^{(\ell)}(\mathbf{t}^{\prime})|^{2}\\
+E\mathbf{|}D\mathbf{_{\mathbf{0}}^{(\ell)}(\mathbf{t}^{\prime})-}%
D\mathbf{_{\mathbf{0}}^{(\ell)}(\mathbf{t})|}^{2}+E\mathbf{|}D_{\mathbf{0}%
}^{(\ell)}(\mathbf{t)-}D_{\mathbf{0}}(\mathbf{t)|}^{2}).
\end{gather*}
Note now that $D_{\mathbf{0}}^{(\ell)}(\mathbf{t}^{\prime})$ is continuous in
$\mathbf{t}^{\prime}$ and so $\lim_{\mathbf{t}^{\prime}\rightarrow\mathbf{t}%
}D_{\mathbf{0}}^{(\ell)}(\mathbf{t}^{\prime})=D\mathbf{_{\mathbf{0}}^{(\ell
)}(\mathbf{t}).}$ Therefore, by taking into account (\ref{limit}), and letting
first $\mathbf{t}^{\prime}\rightarrow\mathbf{t}$ and then $\ell\rightarrow
\infty,$ we obtain for $\lambda^{2}-$almost all ${\mathbf{t\in}}I$ the
approximation%
\[
\lim_{n\rightarrow\infty}\frac{1}{n^{2}}E|S_{\mathbf{n}}(\mathbf{t)}%
-M_{\mathbf{n}}(\mathbf{t)|}^{2}=0.
\]
Furthermore if $\mathbf{n}=(n_{1},n_{2})$ and $f(\mathbf{u})\ln^{+}%
f(\mathbf{u})$ is integrable, by replacing in the proof the limit given in
(\ref{limvar}) by (\ref{limvarln}),\ for $\lambda^{2}-$almost all
${\mathbf{t\in}}I$ we have
\begin{equation}
\lim_{n_{1}>n_{2}\rightarrow\infty}\frac{1}{n_{1}n_{2}}E|S_{\mathbf{n}%
}(\mathbf{t)}-M_{\mathbf{n}}(\mathbf{t)|}^{2}=0. \label{mart approx}%
\end{equation}

\bigskip

By using Theorem 25.4 in Billingsley (1995), the limit (\ref{mart approx})
shows that, the proof of Theorem \ref{CLT} is now reduced to prove the central
limit theorem for $M_{\mathbf{n}}(\mathbf{t).}$

\bigskip

\textbf{The central limit theorem for the martingale.}

\bigskip

\begin{proposition}
\label{propmart}Consider $M_{\mathbf{n}}(\mathbf{t)}$ defined by (\ref{def M})
where $\mathbf{n}=(n_{1},n_{2})$. Then the real and imaginary part of
$M_{\mathbf{n}}(\mathbf{t})$ converge to independent normal random variables
with variance $E|D_{0,0}(\mathbf{t})|^{2}/2$ when $n_{1}\wedge n_{2}%
\rightarrow\infty.$
\end{proposition}

Proof. To ease the notation we shall drop $\mathbf{t}$ and denote
$D_{j,k}=D_{j,k}(\mathbf{t}),$ $M_{\mathbf{n}}=M_{\mathbf{n}}(\mathbf{t}).$

We start by writing ($\mathbf{t=(}t_{1},t_{2})$)
\[
\frac{1}{\sqrt{n_{1}n_{2}}}M_{\mathbf{n}}=\frac{1}{\sqrt{n_{1}}}%
\sum\nolimits_{j=1}^{n_{1}}e^{\mathrm{i}jt_{1}}\frac{1}{\sqrt{n_{2}}}%
\sum\nolimits_{k=1}^{n_{2}}e^{\mathrm{i}kt_{2}}D_{j,k}.
\]
Note that, by construction and since the filtration $(\mathcal{F}_{j}%
,_{k}\mathbf{)}$ is commuting, the sequence $(D_{n_{2},k}^{\prime}%
\mathbf{)}_{k}$ defined by
\[
D_{n_{2},j}^{\prime}\mathbf{=}\frac{1}{\sqrt{n_{2}}}\sum\nolimits_{k=1}%
^{n_{2}}e^{\mathrm{i}kt_{2}}D_{j,k}%
\]
is a triangular array of complex martingale differences with respect to the
filtration $(\mathcal{F}_{j},_{\infty}\mathbf{)}_{j}\mathbf{.}$

For $a$ and $b$ real numbers let us find the limiting distribution of $\ $%
\begin{gather}
a\frac{1}{\sqrt{n_{1}n_{2}}}\operatorname{Re}M_{\mathbf{n}}+b\frac{1}%
{\sqrt{n_{1}n_{2}}}\operatorname{Im}M_{\mathbf{n}}=\frac{1}{\sqrt{n_{1}}}%
\sum\nolimits_{j=1}^{n_{1}}[(a\cos jt_{1}+b\sin jt_{1})\operatorname{Re}%
D_{n_{2},j}^{\prime}\label{sum}\\
+(b\cos jt_{1}-a\sin jt_{1})\operatorname{Im}D_{n_{2},j}^{\prime})]=\frac
{1}{\sqrt{n_{1}}}\sum\nolimits_{j=1}^{n_{1}}\Delta_{n_{2},j}.\nonumber
\end{gather}
In order to find the limiting distribution of $n_{1}^{-1/2}\sum\nolimits_{j=1}%
^{n_{1}}\Delta_{n_{2},j}$ we have to prove now a central limit theorem for the
triangular array of martingale differences $(\Delta_{n_{2},j})_{j\geq1}.$
According to a classical result, which can be found in G\"{a}nssler and
H\"{a}usler (1979), we have to establish that%
\[
\max_{1\leq j\leq n_{1}}\frac{1}{\sqrt{n_{1}}}|\Delta_{n_{2},j}|\rightarrow
^{L_{2}}0\text{ as }n_{1}\wedge n_{2}\rightarrow\infty.
\]
and to verify the Raikov type condition, namely%

\begin{equation}
\frac{1}{n_{1}}E|\sum\nolimits_{j=1}^{n_{1}}(\Delta_{n_{2},j}^{2}%
-E\Delta_{n_{2},j}^{2})|\rightarrow0\text{ as }n_{1}\wedge n_{2}%
\rightarrow\infty. \label{mainlim}%
\end{equation}
The first condition is easy to verify since, by the stationarity involved in
the model and the main result in Peligrad and Wu (2010), the variables
$(|D_{n_{2},j}^{\prime}|^{2})_{j}$ are uniformly integrable and therefore,
$(|\Delta_{n_{2},j}|^{2})_{j}$ in (\ref{sum})\ are also uniformly integrable.

In order to verify (\ref{mainlim}), after using the well-known trigonometric
formulas
\begin{gather*}
2\cos^{2}x=1-\cos2x,\text{ \ \ }2\sin^{2}x=1-\cos2x\\
\cos^{2}x-\sin^{2}x=\cos2x,\text{\ \ \ }2(\cos x)(\sin x)=\sin2x,
\end{gather*}
by Lemma \ref{Lm1} in the Section 6, it follows that for almost all
$\mathbf{t\in}I,$ the terms involving $\cos2jt_{1}$ or $\sin2jt_{1}$ in
(\ref{mainlim}) are negligible as $n_{1}\wedge n_{2}\rightarrow\infty$.

After simple computations, proving (\ref{mainlim}) is reduced to show that
\[
\frac{1}{n_{1}}E|\sum\nolimits_{j=1}^{n_{1}}(|D_{n_{2},j}^{\prime}%
\mathbf{|}^{2}\mathbf{-}E|D_{n_{2},j}^{\prime}|^{2})|\mathbf{\rightarrow}0.
\]
We shall apply Lemma \ref{lm2} below. Clearly it is enough to prove that
\[
\frac{1}{n_{1}}E|\sum\nolimits_{j=1}^{n_{1}}(\operatorname{Re}^{2}D_{n_{2}%
,j}^{\prime}\mathbf{-}E\operatorname{Re}^{2}D_{n_{2},j}^{\prime}%
)|\mathbf{\rightarrow}0\text{ }%
\]
and%
\[
\frac{1}{n_{1}}E|\sum\nolimits_{j=1}^{n_{1}}(\operatorname{Im}^{2}D_{n_{2}%
,j}^{\prime}\mathbf{-}E\operatorname{Im}^{2}D_{n_{2},j}^{\prime}%
)|\mathbf{\rightarrow}0.
\]
Their proofs are similar and we shall deal only with the first one involving
the real part.

Let $m$ be a positive integer. By using Cram\`{e}r theorem, trigonometric
formulas, the main Theorem in Peligrad and Wu (2010), ergodicity of $S$ and
Lemma \ref{Lm1} from Section 6, we can easily show that the vector valued
sequence of martingales $(\operatorname{Re}D_{n_{2},1}^{\prime}%
,...,\operatorname{Re}D_{n_{2},m}^{\prime})_{n_{2}}$ converges to a Gaussian
vector $(N_{1},...,N_{m})$ with the covariance structure
\[
\mathrm{cov}(N_{1},N_{j})=\mathrm{cov}(N_{k},N_{k+j}).
\]
The computations are simple and left to the reader. Because of the martingale
property, we have the orthogonality of $(\operatorname{Re}D_{n_{2},k}^{\prime
})_{k}$. In addition we also have uniform integrability of the family
$(|D_{n_{2},1}^{\prime}|^{2})_{n_{2}}$ , provided by the results in Peligrad
and Wu (2010). By applying the continuous mapping theorem and the convergence
of moments theorem associated to the convergence in distribution (Theorem
25.12 in Billingsley 1995), we obtain
\[
\mathrm{cov}(N_{1},N_{k})=\lim_{n_{2}\rightarrow\infty}\mathrm{cov}%
(\operatorname{Re}D_{n_{2},0}^{\prime},\operatorname{Re}D_{n_{2},k}^{\prime
})=0.
\]
This shows that the Gaussian limit $(N_{k})_{k}$ is a stationary and
independent sequence. It follows that, for all $m\in N,$ $(\operatorname{Re}%
^{2}D_{n_{2},1}^{\prime},...,\operatorname{Re}^{2}D_{n_{2},m}^{\prime}%
)_{n_{2}}$ converges to an independent vector $(N_{1}^{2},...,N_{m}^{2})$ and
$(N_{k}^{2})_{k}$ is stationary and ergodic. Therefore, (\ref{mainlim})\ holds
by Lemma \ref{lm2} in Section 6.

By all of the above considerations we obtain%
\[
a\frac{1}{\sqrt{n_{1}n_{2}}}\operatorname{Re}M_{\mathbf{n}}(\mathbf{t}%
)+b\frac{1}{\sqrt{n_{1}n_{2}}}\operatorname{Im}M_{\mathbf{n}}(\mathbf{t}%
)\Rightarrow(a^{2}+b^{2})N(0,E|D_{0,0}(\mathbf{t})|^{2}),
\]
and the result follows. \ $\square$

\section{Random fields with multi-dimensional index set}

In this section we discuss the differences which occur when the index set is
$Z^{d}.$ The main difference is that we use some recent results on summability
of multi-dimensional trigonometric Fourier series which are surveyed and
further developed in Weisz (2012). Many summability results, needed for our
proofs, have already been extended from dimension $1$ to dimension $d,$ but
the results are very different depending on the summation type and on the
shape of the regions in $Z^{d}$ containing the indexes of summations. Our
intention is to present a method rather than the most general results. The
probabilistic tools are completely developed in our paper. However, the
statements are limited by the level of knowledge in harmonic analysis. In
order to construct the approximating martingale we can always base ourselves
on the summation on cubes, where the celebrated Carleson-Hunt theorem extends
completely for square integrable functions (see Theorem 4.4. in Weisz, 2012)
or we can use an iterative procedure. However, the statements of the CLT\ and
the conditions imposed to the spectral density, strongly depend on shape of
the summation region and the extensions of the Fej\'{e}r-Lebesgue theorem,
namely on the validity of (\ref{limvar}). These regions of summation can be
restricted by using conditions imposed to various norms on $Z^{d}$ or the
summations can be taken over nonrestricted rectangles. In the latter case,
additional restrictions have to be imposed to the spectral density. This is an
active field of research in harmonic analysis and our results can be
reformulated whenever a progress is achieved. We shall formulate the general
results by using only summations over cubes and nonrestricted rectangles.

To introduce the regularity conditions we shall start with a strictly
stationary real valued random field $\mathbf{\xi}=(\xi_{\mathbf{u}%
})_{\mathbf{u}\in Z^{d}}$, defined on the canonical probability space
$R^{Z^{d}}$ and\ define the filtrations $\mathcal{F}_{\mathbf{u}}=\sigma
(\xi_{\mathbf{j}}:\mathbf{j}\leq\mathbf{u})$. Recall that by $\mathbf{j}%
\leq\mathbf{u}$ we understand that each coordinate of $\mathbf{j}$ is less or
equal the corresponding coordinate of $\mathbf{u}.$ By taking intersections of
sigma algebras or sigma algebra generated by unions of sigma algebras, we can
consider the coordinates of $\mathbf{u}$ in $\mathcal{F}_{\mathbf{u}}$ being
valued in $Z\cup\{-\infty,\infty\}.$ The filtration is commuting if
$E_{\mathbf{u}}E_{\mathbf{a}}X=E_{\mathbf{u}\wedge\mathbf{a}}X,$ where the
minimum is taken coordinate-wise. We define $X_{\mathbf{0}}=f((\xi
_{\mathbf{u}})_{\mathbf{u}\in Z^{d}}).$ We also define $T_{i}$ the
coordinate-wise translations and $X_{\mathbf{k}}=f(T_{1}^{k_{1}}\circ...\circ
T_{d}^{k_{d}}(\xi_{\mathbf{u}})_{\mathbf{u}\in Z^{d}}).$ We call $f$ regular
if $E(X_{\mathbf{0}}|\mathcal{F}_{\mathbf{u}})=0$ $a.s.,$ when at least a
coordinate of $\mathbf{u}$ is $-\infty.$

Our general result is summarized in the following theorem:

\begin{theorem}
\label{CLT d}Assume that $\left(  X_{{\mathbf{k}}}\right)  _{{\mathbf{k}}%
\in{Z}^{d}}$ and $(\mathcal{F}_{{\mathbf{k}}})_{{\mathbf{k}}\in{Z}^{d}}$ are
as above and $f$ is regular. In addition, assume that one of the shifts
$T_{i}$ is ergodic,$1\leq i\leq d$. Then, for $\lambda^{d}-$almost all
$\mathbf{t\in}[0,2\pi)^{d},$
\[
\frac{1}{n^{d/2}}(\operatorname{Re}S_{\mathbf{n}}(\mathbf{t}%
),\operatorname{Im}S_{\mathbf{n}}(\mathbf{t}))\Rightarrow(N_{1},N_{2})\text{
as }n\rightarrow\infty,
\]
where $\mathbf{n=}(n,...,n),$ $N_{1},N_{2}$ are i.i.d. normally distributed
random variables with mean $0$ and variance $2^{d-1}\pi^{d}f(\mathbf{t}),$ and
$f(\mathbf{t})$ is the spectral density of the sequence $\left(
X_{{\mathbf{k}}}\right)  _{{\mathbf{k}}\in{Z}^{d}}$. Furthermore, if
$f(\mathbf{u})(\ln^{+}f(\mathbf{u}))^{d-1}$ is integrable, then%
\[
\frac{1}{\sqrt{n_{1}n_{2}...n_{d}}}(\operatorname{Re}S_{\mathbf{n}}%
(\mathbf{t}),\operatorname{Im}S_{\mathbf{n}}(\mathbf{t}))\Rightarrow
(N_{1},N_{2})\text{ as }\wedge_{1\leq i\leq d}n_{i}\rightarrow\infty,
\]
with $N_{1},N_{2}$ as above.
\end{theorem}

\textbf{Proof of Theorem \ref{CLT d}. }The proof of this theorem follows the
same lines as of Theorem \ref{CLT} with the following differences. In order to
be able to obtain a characterization of the spectral density, we have to
introduce the $d$-dimensional projection operator. By using the commutative
property of the filtrations it is convenient to define%

\[
\mathcal{P}_{\mathbf{0}}(X)=\mathcal{P}_{1}\circ\mathcal{P}_{2}\circ
...\circ\mathcal{P}_{d}(X),
\]
where%
\[
\mathcal{P}_{j}(Y)=E(Y|\mathcal{F}_{0}^{(j)})-E(Y|\mathcal{F}_{-1}^{(j)}).
\]
Above we used the notation: $\mathcal{F}_{0}^{(j)}=\mathcal{F}_{\mathbf{0}}$,
and $\mathcal{F}_{-1}^{(j)}=\mathcal{F}_{\mathbf{u}}$, where $\mathbf{u}$ has
all the coordinates $0$ with the exception of the $j$-th coordinate, which is
$-1$. For instance when $d=3,$ $\mathcal{P}_{2}(Y)=E(Y|\mathcal{F}%
_{0,0,0})-E(Y|\mathcal{F}_{0,-1,0}).$

We can easily see that, by using the commutativity property, this definition
is a generalization of the case $d=2$. We note that, by using this definition
of $\mathcal{P}_{\mathbf{0}}(X),$ the statement and the proof of Theorem
\ref{specden} remain unchanged if we replace $Z^{2}$ with $Z^{d}.$ The
definition of the approximating martingale is also clear as well as the proof.
We point out the following two differences in the proof. One difference is
that instead of Theorem 1 in Marcinkiewicz and Zygmund (1939) we use Corollary
14.4 in Weisz (2012), which assures the validity of\ (\ref{limvar}) for
$\lambda^{d}-$almost all $\mathbf{t}$ in $[0,2\pi)^{d}.$ Another difference in
the proof is that Proposition \ref{propmart} is proved by induction. More
precisely, we use instead of the results in Peligrad and Wu (2010) the
induction hypothesis. For proving the second part of the Theorem \ref{CLT d},
we use several results in Weisz (2012), namely Corollary 16.5 about
unrestricted summability and the line above relation 15.2 on page 123.

\section{Examples}

We start this section by mentioning an easy way to generate commuting
filtrations. This happens for instance when we consider a stationary random
field $\mathbf{\xi=(\bar{\xi}}_{\ell})_{\ell\in Z}$ with its columns
$\mathbf{\bar{\xi}}_{\ell}=(\mathbf{\xi}_{u,\ell})_{u\in Z}$ independent
copies of a stationary stochastic process. Of course, we can also consider as
well random fields with independent lines. This can be seen by combining the
Lemmas (\ref{indeprows}) and (\ref{commuting}) in Section 6. The random field
of interest is then constructed by taking functions of $(\xi_{k,j})_{k,j\in
Z}$ as in definition (\ref{defXfield}). Furthermore, if the columns of
$(\xi_{k,j})$ are independent, then $\mathcal{F}_{0,-\infty}$ is trivial. If
the lines of $(\xi_{k,j})_{k,j\in Z}$ are independent then $\mathcal{F}%
_{-\infty,0}$ is trivial.

\bigskip

Next, we give examples of stationary random fields $(\xi_{\mathbf{u}%
})_{\mathbf{u}\in Z^{2}}$ which generate commuting filtration and in addition
both $\mathcal{F}_{0,-\infty}$ and $\mathcal{F}_{-\infty,0}$ are trivial.

\bigskip

\textbf{Independent copies of a stationary sequence with "nonparallel" past
and future. }

\bigskip

The $\rho-$mixing coefficient, also known as maximal coefficient of
correlation is defined as
\[
\rho(\mathcal{A},\mathcal{B})=\sup\{\mathrm{Cov}(X,Y)/\Vert X\Vert_{2}\Vert
Y\Vert_{2}:\,X\in L_{2}(\mathcal{A}),\,Y\in L_{2}(\mathcal{B})\}.
\]
For the stationary sequence of random variables $(\xi_{k})_{k\in Z}$, denote
by $\mathcal{F}_{0}$ the past $\sigma$--field generated by $\xi_{k}$ with
indices $k\leq0$ and by $\mathcal{F}^{n}$ the future $\sigma$--field after
$n-$steps generated by $\xi_{j}$ with indices $j\geq n.$ The sequence of
coefficients $(\rho_{n})_{n\geq1}$ is then defined by
\[
\rho_{n}=\rho(\mathcal{F}_{0},\mathcal{F}^{n}\mathcal{)}.
\]
If $\rho_{n}<1$ for some $n>1,$ then the tail sigma field $\mathcal{F}%
_{-\infty}=\cap_{n\in Z}\mathcal{\sigma((}\xi_{j})_{j\leq n})$ is trivial; see
Proposition (5.6) in Bradley (2007). In this case it is customary to say that
$\mathcal{F}_{0}$ and $\mathcal{F}^{n}$ are not parallel. Now we take a random
field with columns $\bar{\xi}_{j}=(\xi_{k,j})_{k\in Z}$ independent copies of
a stationary sequence with $\rho_{n}<1$ for some $n>1$. Clearly, because the
columns are independent the sigma field $\mathcal{F}_{0,-\infty}$ is trivial.
Furthermore, we shall argue that $\mathcal{F}_{-\infty,0}$ is also trivial. To
prove it, we apply Theorem 6.2 in in Csaki and Fisher (1963) (see also Theorem
6.1 in Bradley, 2007). According to this theorem%
\begin{align*}
&  \rho(\sigma((\xi_{k,j})_{j\in Z,k\leq0}),\sigma((\xi_{k,j})_{j\in Z,k\geq
n})\\
&  =\sup_{j}\rho(\sigma((\xi_{k,j})_{k\leq0}),\sigma((\xi_{k,j})_{k\geq
n}))=\rho_{n}<1.
\end{align*}
Therefore we also have that $\mathcal{F}_{-\infty,0}$ is trivial and our
theorem applies. For this case we obtain the following corollary:

\begin{corollary}
Assume that the random field $(\xi_{k,j})_{k,j\in Z}$ consists of columns,
which are independent copies of a stationary sequence $(\xi_{j})_{j\in Z}$
having $\rho_{n}<1$ for some $n>1.$ Construct $(X_{\mathbf{n}})_{\mathbf{n\in
Z}^{2}}$ by (\ref{defXfield}) and assume that the variables are centered and
square integrable. Then the results of Theorems \ref{specden} and \ref{CLT} hold.
\end{corollary}

As a particular example we can take, as generator of the commuting sigma
algebras, independent copies of a Gaussian sequence with a special type of
spectral density. It is convenient to define the spectral density on the unit
circle in the complex plane, denoted by $T$. Let $%
\mu
$ denote normalized Lebesgue measure on $T$ (normalized so that $%
\mu
(T)=1$). For a given random sequence $X:=(X_{k})_{k\in Z}$, a
\textquotedblleft spectral density function\textquotedblright\ (if one exists)
can also be viewed as a real, nonnegative, Borel, integrable function $f$
$:T\rightarrow\lbrack0,\infty)$ such that for every $k\in Z$%

\[
\mathrm{cov}(X_{k},X_{0})=\int\limits_{t\in T}t^{k}f(t)%
\mu
(dt).
\]

Let $(\xi_{j})_{j\in Z}$ be a stationary Gaussian sequence and let $n$ be a
positive integer. The following two conditions are equivalent:

(a) $\rho_{n}<1$.

(b) $(\xi_{j})_{j\in Z}$ has a spectral density function $f$ (on $T$) of the form%

\[
f(t)=|p(t)|exp(u(t)+\tilde{v}(t)),\text{ }t\in T
\]
where $p$ is a polynomial of degree at most $n-1$ (constant if $n=1$), $u$ and
$v$ are real bounded Borel functions on $T$ with $||v||_{\infty}<\pi/2$, and
$\tilde{v}$ is the conjugate function of $v$.

For $n=1$, this equivalence is due to Helson and Szeg\"{o} (1960). For general
$n\geq1,$ it is due to Helson and Sarason (Theorem 6, 1967).

\bigskip

\textbf{Functions of i.i.d.}

\bigskip

Our results also hold for any random field which is Bernoulli, i.e. a function
of i.i.d. random field. For instance, if $(\xi_{\mathbf{n}})_{\mathbf{n}\in
Z^{d}}$ is a random field of independent, identically distributed random
variables and we define $(X_{\mathbf{k}})_{\mathbf{k\in Z}^{d}}$ and
$(\mathcal{F}_{{\mathbf{k}}})_{{\mathbf{k}}\in{Z}^{d}}$ as in Theorem
\ref{CLT d}.\ Then the filtration is commuting and the regularity conditions
of Theorem \ref{CLT d} are satisfied provided the variables are centered. If
in addition $X_{\mathbf{0}}$ is square integrable, then the result of Theorem
\ref{CLT d} holds.

For the next two examples the only conditions imposed are equivalent to the
existence of the fields involved.

\bigskip

\begin{example}
\label{exinear}(Linear field) Let $(\xi_{\mathbf{n}})_{\mathbf{n}\in Z^{d}}$
be a random field of independent, identically distributed random variables
which are centered and have finite second moment. Define%
\[
X_{\mathbf{k}}=\sum_{\mathbf{j}\in Z^{d}}a_{\mathbf{k}-\mathbf{j}}%
\xi_{\mathbf{j}}.
\]
Assume that $\sum_{\mathbf{j}\in Z^{d}}a_{\mathbf{j}}^{2}<\infty$. Then the
CLT in Theorem \ref{CLT d} holds.
\end{example}

\bigskip

Another class of nonlinear random fields is the Volterra process, which plays
an important role in the nonlinear system theory.

\begin{example}
\label{Volterra}(Volterra field) Let $(\xi_{\mathbf{n}})_{\mathbf{n}\in Z^{d}%
}$ be a random field of independent random variables identically distributed
centered and with finite second moment. Define%
\[
X_{\mathbf{k}}=\sum_{\mathbf{u},\mathbf{v}\in Z^{d}}a_{\mathbf{u},\mathbf{v}%
}\xi_{\mathbf{k-u}}\xi_{\mathbf{k-v}}\text{ },
\]
where $a_{\mathbf{u},\mathbf{v}}$ are real coefficients with $a_{\mathbf{u}%
,\mathbf{u}}=0$ and $\sum_{\mathbf{u},\mathbf{v}\in Z^{d}}a_{\mathbf{u,v}}%
^{2}<\infty.$ Then the CLT in Theorem \ref{CLT d} holds.
\end{example}

\section{Supplementary results}

In this section we prove two auxiliary results. They are laws of large numbers
which have interest in themselves.

The following lemma is an extension of a result in Zhang (2017).

\bigskip

\begin{lemma}
\label{Lm1}Assume that the triangular array $(X_{n_{2},k})_{k\in Z}$ is
row-wise stationary and $(X_{n_{2},k})_{n_{2}\geq1}$ is uniformly integrable
for any $k$ fixed. In addition assume that $X_{n_{2},k}\Rightarrow X_{k}$,
where $X_{k}$'s have the same distribution and are in $L_{1}$. Then, for
$\lambda-$almost all $t\in\lbrack-\pi,\pi)$,%
\[
\frac{1}{n_{1}}\sum\nolimits_{k=1}^{n_{1}}e^{\mathrm{i}kt}X_{n_{2}%
,k}\rightarrow0\text{ }a.s.\text{ when }n_{1}\wedge n_{2}\rightarrow\infty.
\]

\end{lemma}

Proof. Let $m\geq1$ be a fixed integer and define consecutive blocks of
indexes of size $m$, $I_{j}(m)=\{(j-1)m+1,...,mj\}.$ In the set of integers
from $1$ to $n$ we have $k_{n_{1}}=k_{n_{1}}(m)=[n_{1}/m]$ such blocks of
integers and a last one containing less than $m$ indexes. By the uniform
integrability of $(X_{n_{2},k})_{n_{2}\geq1}$ we have%
\begin{align*}
\frac{1}{n_{1}}E|\sum\nolimits_{k=1}^{n_{1}}e^{\mathrm{i}kt}X_{n_{2},k}|  &
\leq\frac{1}{n_{1}}\sum\nolimits_{j=1}^{k_{n_{1}}(m)}E|\sum\nolimits_{k\in
I_{j}(m)}e^{\mathrm{i}kt}X_{n_{2},k}|\\
&  +\frac{1}{n_{1}}E|\sum\nolimits_{j=k_{n_{1}}(m)+1}^{n_{1}}e^{\mathrm{i}%
jt}X_{n_{2},j}|\\
&  \leq\frac{1}{m}E|\sum\nolimits_{k=1}^{m}e^{\mathrm{i}kt}X_{n_{2}%
,k}|+o_{n_{1}}(1)\text{ as }n_{1}\rightarrow\infty.
\end{align*}
Now, again by the uniform integrability of $(X_{n_{2},k})_{n_{2}\geq1}$ and
the convergence of moments associated to the weak convergence (see Theorem
25.12 in Billingsley, 1995) we have
\begin{gather*}
\lim_{n_{2}\rightarrow\infty}E|\sum\nolimits_{k=1}^{m}e^{\mathrm{i}kt}%
X_{n_{2},k}|\leq\\
\lim_{n_{2}\rightarrow\infty}E|\sum\nolimits_{k=1}^{m}X_{n_{2},k}\sin
kt|+\lim_{n_{2}\rightarrow\infty}E|\sum\nolimits_{k=1}^{m}X_{n_{2},k}\cos
kt|\\
=E|\sum\nolimits_{k=1}^{m}X_{k}\sin kt|+E|\sum\nolimits_{k=1}^{m}X_{k}\cos
kt|.
\end{gather*}
Since the $X_{k}$'s have the same distribution, by Zhang (2017), for almost
all $t\in\lbrack-\pi,\pi)$
\[
\frac{1}{m}\sum\nolimits_{k=1}^{m}X_{k}\sin kt\rightarrow0\text{ and }\frac
{1}{m}\sum\nolimits_{k=1}^{m}X_{k}\cos kt\rightarrow0\text{ }P-a.s.\text{ and
in }L_{1}.
\]
$\square$

\bigskip

\begin{lemma}
\label{lm2}Assume that the triangular array $(X_{n_{2},k})_{k\in Z}$ is
row-wise stationary, mean $0$ and $(X_{n_{2},k})_{n_{2}\geq1}$ is uniformly
integrable for any $k$ fixed. In addition assume that the finite dimensional
distributions of $(X_{n_{2},k})_{k}$ converge in distribution to those of
$(X_{k})_{k}$ as $n_{2}\rightarrow\infty,$ where $(X_{k})_{k}$ is stationary
and ergodic and in $L_{1}.$ Then
\[
\frac{1}{n_{1}}\sum\nolimits_{k=1}^{n_{1}}X_{n_{2},k}\ \text{converges in
}L_{1}\text{ to }0\text{ when }n_{1}\wedge n_{2}\rightarrow\infty.
\]

\end{lemma}

Proof. As in the previous lemma, we make blocks of variables as before and use
the inequality
\[
\frac{1}{n_{1}}E|\sum\nolimits_{k=1}^{n_{1}}X_{n_{2},k}|\leq\frac{1}{m}%
E|\sum\nolimits_{k=1}^{m}X_{n_{2},k}|+o_{n_{1}}(1)\text{ as }n_{1}%
\rightarrow\infty.
\]
Now, by using the uniform integrability of $(X_{n_{2},k})_{n\geq1}$ we obtain
\[
\lim_{n\rightarrow\infty}E|\sum\nolimits_{k=1}^{m}X_{n_{2},k}|=E|\sum
\nolimits_{k=1}^{m}X_{k}|.
\]
Furthermore, note also that by the conditions of this lemma we also have
$E(X_{k})=0$ for all $k.$ By the ergodic theorem $\sum\nolimits_{k=1}^{m}%
X_{k}/m\rightarrow0$ $a.s.$ and in $L_{1}$ and therefore%
\[
E|\frac{1}{m}\sum\nolimits_{k=1}^{m}X_{k}|\rightarrow0\text{ as }%
m\rightarrow\infty.
\]
$\square$

\bigskip

The following lemma follows by Problem 34.11 in Billingsley (1995).

\bigskip

\begin{lemma}
\label{commuting}Assume that $X,Y,Z$ are integrable random variables. Then the
following are equivalent%
\[
E(g(X,Y)|\sigma(Y,Z))=E(g(X,Y)|Y)\text{ }a.s.
\]%
\[
E(g(Z,Y)|\sigma(Y,X))=E(g(Z,Y)|Y)\text{ }a.s.
\]

\end{lemma}

We recall the following lemma which is not difficult to verify. Its proof is
left to the reader.{}

\bigskip

\begin{lemma}
\label{indeprows}Assume that $X,Y,Z$ are integrable random variables such that
$(X,Y)$ and $Z$ are independent. Assume that $g(X,Y)$ is integrable. Then%
\[
E(g(X,Y)|\sigma(Y,Z))=E(g(X,Y)|Y)\text{ }a.s.
\]

\end{lemma}

\textbf{Acknowledgements.} This research was supported in part by the
NSF\ grant DMS-1512936. The authors are grateful to the referees for numerous
suggestions, which contributed to a significant improvement of a previous
version of the paper.


\begin{thebibliography}{99}                                                                                               %


\bibitem {b}Bary, N. K. (1964). A Treatise on Trigonometric Series\textit{.}
New York, Macmillan.

\bibitem {bi}Billingsley, P. (1995). Probability and measures. (3rd ed.).
Wiley Series in Probability and Statistics, New York.

\bibitem {Br}Bradley, R.C. (2007). Introduction to strong mixing conditions. 3
Volumes, Kendrick Press.

\bibitem {Ca}Carleson, L. (1966). On convergence and growth of partial sums of
Fourier series. Acta Math\textit{.} 116 135-157.

\bibitem {CC}Cohen, G. and J.P. Conze (2013). The CLT for rotated ergodic sums
and related processes. Discrete and Continuous Dynamical Systems - Series A
(DCDS-A) 33 3981 - 4002.

\bibitem {CF}Cs\'{a}ki P. and J. Fischer (1963). On the general notion of
maximal correlation. Magyar Tud. Akad. Mat. Kutat\'{o} Int. K\"{o}zl. 8 27-51.

\bibitem {CDV}Cuny, C., Dedecker, J. and D. Voln\'{y} (2015). A functional
central limit theorem for fields of commuting transformations via martingale
approximation. Zapiski Nauchnyh Seminarov POMI 441.C.and J. Math. Sci. (N.Y.)
2016, 219, (5) 765--781.

\bibitem {EM}El Machkouri, M., Voln\'{y}, D. and W.B. Wu (2013). A central
limit theorem for stationary random fields. Stochastic Process. Appl. 123 1-14.

\bibitem {cf1}Fefferman, C. (1971 a). On the convergence of multiple Fourier
series. Bull. Amer. Math. Soc. 77 744--745.

\bibitem {cf2}Fefferman, C. (1971 b). On the divergence of multiple Fourier
series. Bull. Amer. Math. Soc. 77 191--195.

\bibitem {Franco}Francos, J.M., Meiri, A.Z. and B. Porat (1995). A Wold-like
decomposition of two-dimensional discrete homogeneous random fields. Ann.
Appl. Probab. 5 248-260.

\bibitem {GH}G\"{a}nssler, P. and E. H\"{a}usler (1979). Remarks on the
functional central limit theorem for martingales. Z.
Wahrscheinlichkeitstheorie verw. Gebiete 50 237-243.

\bibitem {G}Giraudo, D. (2017). Invariance principle via orthomartingale
approximation. arXiv:1702.08288.

\bibitem {G}Gordin, M. I. (2009). Martingale-coboundary representation for a
class of stationary random fields. Zap. Nauchn. Sem. S.-Peterburg. Otdel. Mat.
Inst. Steklov. (POMI) 364, Veroyatnostn i Statistika. 14.2, 88-108, 236; and
J. Math. Sci. (N. Y.) 163 (2009) 363-374.

\bibitem {H2}Helson, H. and G. Szeg\"{o}. (1960). A problem in prediction
theory. Ann. Mat. Pura Appl. 51 107-138.

\bibitem {Helson}Helson, H. and D. Lowdenslager (1962). Prediction theory and
Fourier series in several variables. Acta Math. 106 175-213.

\bibitem {H3}Helson, H. and D. Sarason. (1967). Past and future. Math. Scand.
21 5-16.

\bibitem {hunt}Hunt, R. I. and Young, W. S. (1974). A weighted norm inequality
for Fourier series, Bull. Amer. Math. Soc. 80, 274--277.

\bibitem {J}Jessen, B., Marcinkiewicz, J. and A. Zygmund (1935). Note on the
differentiability of multiple integrals. Fundamenta Mathematicae 25 217-234.

\bibitem {KaM}Kallianpur, G., Miamee, A.G. and H. Niemi (1990). On the
prediction theory of two-parameter stationary random fields. J. Multiv. Anal.
32 120-149.

\bibitem {lp}Lifshits, M. A.and M. Peligrad (2015). On the spectral density of
stationary processes and random fields. Zap. Nauchn. Sem. S.-Peterburg. Otdel.
Mat. Inst. Steklov. (POMI) 441, Veroyatnost' i Statistika. 22, 274--285 and J.
Math. Sci. (N.Y.) 219 (2016) 789--797.

\bibitem {MZ}Marcinkiewicz, J. and A. Zygmund (1939). On the summability of
double Fourier series. Fundamenta Mathematicae 32 122-132.

\bibitem {M}Miller, C. (1995). A CLT\ for the periodograms of a $\rho^{\ast}%
-$mixing random field. Stochastic Processes and their Applications 60 313--330.

\bibitem {PeWu}Peligrad, M. and W. B. Wu (2010). Central limit theorem for
Fourier transforms of stationary processes. Ann. Probab\textit{.} 38 2009-2022.

\bibitem {PZ}Peligrad, M. and Na Zhang (2017). On the normal approximation for
random fields via martingale methods. arXiv:1702.01143.

\bibitem {R}Rosenblatt, M. (1972). Central limit theorem for stationary
processes. Berkeley Symp. on Math. Statist. and Prob. Proc. Sixth Berkeley
Symp. on Math. Statist. and Prob., Vol. 2 (Univ. of Calif. Press) 551-561.

\bibitem {Shu}Schuster, A. (1898). On the investigation of hidden
periodicities with application to a supposed 26 day period of meteorological
phenomena. Terrestrial Magnetism and Atmospheric Electricity 3 13-41.

\bibitem {V}Voln\'{y}, D. (2015). A central limit theorem for fields of
martingale differences, C. R. Math. Acad. Sci. Paris 353 1159-1163.

\bibitem {VW}Voln\'{y} D. and Y. Wang (2014). An invariance principle for
stationary random fields under Hannan's condition. Stochastic Process. Appl,
124 4012-4029.

\bibitem {WW}Wang, Y. and M. Woodroofe (2013). A new criteria for the
invariance principle for stationary random fields. Statistica Sinica 23 1673-1696.

\bibitem {Weisz}Weisz, F. (2012). Summability of multi-dimensional
trigonometric Fourier series. Surveys in Approximation Theory 7 1--179.

\bibitem {ZNa}Zhang, Na (2017). On the law of large numbers for discrete
Fourier transform. Statistics \& Probability Letters 120 (C) 101-107.
\end{thebibliography}
\end{document}